\documentclass[12pt,draftcls,onecolumn]{IEEEtran}

\usepackage{graphicx}

\usepackage[cmex10]{amsmath}
\usepackage{amsfonts}
\usepackage{amssymb}

\usepackage{algorithmic}

\usepackage{url}

\usepackage{ulem}

\begin{document}

\title{Direct search methods for an open problem of optimization in systems and control}

\author{Emile~Simon
        ~and~Vincent~Wertz
\thanks{E. Simon and V. Wertz are with the Mathematical Engineering Department, ICTEAM Institute, Universit\'e Catholique de Louvain, 4 avenue Georges Lema\^{i}tre, 1348 Louvain-la-Neuve, Belgium (Tel. +32 10 47 21 80; fax +32 10 47 80 32. e-mail: Emile.Simon@uclouvain.be).}}
\maketitle

\begin{abstract}

The motivation of this work is to illustrate the efficiency of some often overlooked alternatives to deal with optimization problems in systems and control. In particular, we will consider a problem for which an iterative linear matrix inequality algorithm (ILMI) has been proposed recently. As it often happens, this algorithm does not have guaranteed global convergence and therefore many methods may perform better. We will put forward how some general purpose optimization solvers are more suited than the ILMI. This is illustrated with the considered problem and example, but the general observations remain valid for many similar situations in the literature.

\end{abstract}

\begin{IEEEkeywords}
Optimization; Linear systems; Positive filtering
\end{IEEEkeywords}

\section{Introduction}

\subsection{Motivation}
There are many open problems of optimization in systems and control theory, often non-smooth, non-convex and NP-hard or of unknown complexity. Because of the success obtained with linear matrix inequalities (LMIs) since the mid-nineties, a general tendency in the system and control literature is to formulate many of these problems with bilinear matrix inequalities (BMIs) or LMIs plus a non-convex rank constraint. Once these formulations obtained, the usual techniques for trying to solve these problems are iterative linear matrix inequalities algorithms (ILMIs). Such algorithms solve successively different LMI subproblems of the original BMI formulation. However, these algorithms loose the guaranteed polynomial time complexity of LMI solvers but more importantly most do not have guaranteed convergence to locally optimal solutions (i.e. global convergence when starting from any feasible initial solution).

Such ILMIS, heuristics without global convergence, are then likely to be outperformed by other methods like general purpose optimization solvers. This is the fact that we aim at putting forward with the current paper: much research work sticks to a unique approach, like BMIs/LMIs, and ignores alternative formulations which can however be much more appropriate. 

The performance can be compared in three directions: 1) Objective value, 2) Computational time and 3) User time. More details will be given in the paper but we already outline here some key elements. 1) The methods that will be proposed are often found in the literature to behave well on non-smooth functions and some are proved globally convergent on smooth functions. This is not the case for many ILMIs. 2) On non-smooth and non-convex problems there are no worst-case complexity bounds neither for the proposed methods nor for ILMIs. Thus, on this front, only experience will tell. 3) The time required by the user to implement the method is clearly smaller with general purpose optimization solvers, which are intended to minimize any user-defined function $f(\theta):\mathbb{R}^n\rightarrow \mathbb{R}$. Formulating a problem under that form may be done in a matter of minutes or hours. But developing a theoretically consistent ILMI algorithm, along with a working implementation, can easily require several weeks. This is by the way admissible in fundamental research but clearly unacceptable in the industry.

\subsection{Considered problem}
To illustrate the point of the paper, we will consider a problem for which an ILMI approach was recently proposed \cite{LLS10}. This problem is the design of a reduced-order positive filter for linear systems. Positive systems, found in many areas (see \cite{FR00}), are dynamic systems with state variables and outputs positive at all times. The filter to be designed has to estimate an unmeasured output ($z$) of a system ($\Sigma$) from the measurements ($y$), with minimum objective value chosen as the $\mathcal{H}_\infty$ norm of the filtering error transfer function ($\mathcal{G}$) between the exogenous disturbance signal ($w$) and the error ($e$): the difference between the output ($z$) and its estimate ($\hat{z}$). Details on the problem motivations and challenges are given in the introduction of \cite{LLS10}, here only the elements necessary to reproduce the contribution are recalled.

Consider the following asymptotically stable discrete-time system:
\begin{equation}
\Sigma:
\left\{\begin{array}{ccccc}
   x_{k+1} & = & A x_k & + & B w_k  \\
   y_k & = & C x_k & + & D w_k \\
   z_k & = & L x_k & + & G w_k
\end{array}\right.
\label{Sys}
\end{equation}
The notations are classical and identical to those of \cite{LLS10}, where all terms are more formally detailed. This system is positive if and only if the matrices $A,B,C,D,L,G$ have only positive entries \cite{FR00}. Note also that in this paper we only require the positivity of the filter to be designed and not of the system $\Sigma$, but it can be assumed because otherwise it is not necessary to design a positive filter.

The aim of the filter is to compute an estimation $\hat{z}_k$ of the unmeasured signal $z_k$ in $\Sigma$ from the measured signal $y_k$. More specifically, we want to build the following filter:
\begin{equation}
\hat{\Sigma}:
\left\{\begin{array}{ccccc}
   \hat{x}_{k+1} & = & \hat{A} \hat{x}_k & + & \hat{B} y_k  \\
   \hat{z}_k & = & \hat{C} \hat{x}_k & + & \hat{D} y_k
\end{array}\right.
\label{Filt}
\end{equation}
where $\hat{A},\hat{B},\hat{C},\hat{D}$ are the filtering parameters to be determined, matrices having only positive entries to ensure the positiveness of $\hat{z}_k$ \cite{FR00,LLS10}. The difficulty of designing this filter stems from that particular requirement, which prevents the application of methods normally used for this kind of problem (see \cite{LLS10} and ref. therein).

Defining $\xi_k=[x_k^T,\hat{x}_k^T]^T$ and $e_k=z_k-\hat{z}_k$, we get from (\ref{Sys}) and (\ref{Filt}) the description of the filtering error system:
\begin{equation}
\Sigma_e:
\left\{\begin{array}{ccccc}
   \xi_{k+1} & = & A_f \xi_k & + & B_f w_k  \\
   e_k & = & C_f \xi_k & + & D_f w_k
   \end{array}\right.
\label{ErSy}
\end{equation}
where
	\[
	\begin{array}{ll}
    A_f = \begin{bmatrix}
		A & 0 \\
		\hat{B}C & \hat{A}
		\end{bmatrix}, & B_f =  \begin{bmatrix}
		B\\
		\hat{B}D
		\end{bmatrix}
		\\
    C_f  = 
			\begin{bmatrix}
			L-\hat{D}C & -\hat{C}
			\end{bmatrix}, & D_f = G - \hat{D}D
\end{array}
\]
The transfer function of the filtering error system $\Sigma_e$ is given by:
	\[\mathcal{G}(z) = C_f(zI-A_f)^{-1}B_f+D_f
\]
The problem considered is defined hereunder. 

\textit{Reduced-Order Minimal $\mathcal{H}_\infty$ Positive Filtering Problem:} 
\begin{equation}
\min_{\hat{\Sigma}} \ ||\mathcal{G}||_\infty\ s.t.\ \Sigma_e\ stable\ and\ \hat{A}, \hat{B}, \hat{C}, \hat{D}\ positive 
\label{ProbEq}
\end{equation}
This problem is a generalization of the one in \cite{LLS10}. The difference is that we do not only seek a filter respecting a given level $||\mathcal{G}||_\infty <\gamma$ but we rather wish to minimize this level. This objective is more meaningful and also more practical for the user who does not need to choose an arbitrary level $\gamma$ beforehand. Note that considering the expression of $A_f$, since $\Sigma$ is stable we have that the filtering error system $\Sigma_e$ is stable iff the filter $\hat{\Sigma}$ is stable.

\section{The general approach}

The problem above can be expressed as either of the two following optimization problems of minimizing an objective function $f(\theta):\mathbb{R}^n\rightarrow \mathbb{R}$:\\
\textit{The unconstrained problem:}                
\[\min_{\theta}\ f(\theta)\ =\
	\left\{\begin{array}{ll}
	   \infty & if\ \theta\ not\ positive\ or\ \Sigma_e\ unstable  \\
	   ||\mathcal{G}||_\infty & otherwise
	   \end{array}\right.
\]
\textit{The constrained problem:}
	\[\min_{\theta}\ f(\theta)\ =\ ||\mathcal{G}||_\infty\ (\infty\ if\ \Sigma_e\ unstable)\ s.t.\ \theta\ positive
\]
where $n$ is the number of variables in $\theta$, the vector containing the entries of $\hat{A}, \hat{B}, \hat{C}, \hat{D}$. These problems can actually be dealt with many general-purpose optimization solvers. Amongst these we outline some Direct Search (DS) methods, using the unconstrained formulation: the Mesh Adaptive Direct Search (MADS) \cite{AD06}, the Multidirectional Search (MDS) \cite{T89,MCT}, the alternative directions \cite{MCT}, the Nelder-Mead algorithm (NM) \cite{NM65,MCT,LRWW98,PCB02,LL04,BPT06}, the Particle Swarm Optimizer (PSO) \cite{VE02} (see \cite{KLT03} for the latest survey on DS methods). We also mention gradient-based methods, which standard and widely known principles are implemented e.g. in the following two functions of Matlab: \texttt{fminunc} for the unconstrained formulation and \texttt{fmincon} for the constrained formulation, both methods estimating numerical gradients when the expression of the gradient of $f$ is not given (see the Matlab help \cite{Opti} and references therein for detailed descriptions). All these methods can be used for the problem discussed here. 

Note that when the expression of an (estimate) of the gradient is known, more evolved techniques can be developed (like \cite{AN06,GHMO09}) to (try to) reduce the computational time, especially for problems with larger number of variables. For complex problems like the one considered, the gradient information is useful to accelerate the search but without guarantee to find better optima than methods which do not use this information (see benchmarks in \cite{S11}).

\subsection{A method put forward}

For illustration purposes we put forward one method in particular, the Nelder-Mead algorithm (NM). This method is probably the most popular algorithm for unconstrained minimization of problems of moderate dimension (e.g. $<$ 50 to 100 variables). The central question is whether locally optimal solutions can be reached, which is not guaranteed for complex problems with most (if not all) iterative Linear Matrix Inequality algorithms (ILMIs) including the one proposed in \cite{LLS10}, as motivated in the introduction. Note that for further illustration, we will also use several other methods than NM in the next section.

The NM algorithm does not appear often in (recent) systems and control literature. This algorithm was first proposed almost fifty years ago in \cite{NM65} and belongs to the class of DS methods. These methods, characterized by the fact that they do not build a model of the objective function, belong to the broader class of derivative-free optimization methods which do not use any gradient or Hessian information. The reader is referred to the last survey on DS methods \cite[Sec. 1.4]{KLT03} for a broader description. The basic ideas behind NM are briefly described in the next paragraph, detailed descriptions can be found in \cite{LRWW98} or \cite{NM65,PCB02}.

The first step is the generation of an initial simplex of $n+1$ solutions around and including the provided initial solution. The objective function is evaluated at each of these $n+1$ solutions and sorted from the best to the worst. Then NM chooses iteratively between several possibilities (or steps) to change the shape of the simplex (eventually displacing it), trying to find better solutions. For example the basic step is that the worst solution is reflected on the other side of the simplex, in order to create a `downhill' effect.

Two easily available implementations of NM are the \texttt{fminsearch} function in the Optimization Toolbox of Matlab \cite{LRWW98,Opti} or the \texttt{nmsmax} function of \cite{MCT}. It must be noted however that these implementations may fail to converge to locally optimal solutions, starting from a feasible initial solution. Indeed the original NM can fail because of the deterioration of the simplex geometry or lack of sufficient decrease. However the convergence of the method can be guaranteed on smooth functions by taking care of these situations, which is done for example in \cite{PCB02} where an additional step ensuring the convergence is proposed. Actually there is a third easily available implementation, the grid restrained NM algorithm \cite{BPT06}. This algorithm enjoys, like \cite{PCB02}, guaranteed convergence on smooth functions and will also be used in our tests.

An earlier method to improve the convergence of NM proposed in \cite{LL04} consists of an involved restarting strategy. Here we suggest a method with easier implementation: restarting the algorithm until the latest obtained objective value is not better than the previous one to a given accuracy. This does not formally guarantee the convergence to local solutions on smooth objective functions but renders it much more likely. Indeed, restarting NM regenerates its simplex and this allows in practice to properly span the space $\mathbb{R}^n$ around the last solution. 

Anyway, the hypothesis for convergence analysis in \cite{PCB02,LL04,BPT06} require that the function be smooth. And as pointed out in \cite{AN06}, the $\mathcal{H}_\infty$ norm objective function is non-smooth and this may cause the algorithm -or any of the techniques cited above- to stop at suboptimal solutions. This is the reason why we will use not only \cite{BPT06} but also try the idea of restarting NM in the results section.

We summarize that the point is not to guarantee the convergence to locally optimal solutions for all optimizations, which require smooth objective functions, but instead to propose an easily implementable technique that performs well and in practice often reaches locally optimal solutions even for non-smooth objectives functions. 

Nelder-Mead with local restart(s) can be more formally defined as:

\begin{algorithmic}
\STATE $\theta_f(1) = $\text{Nelder-Mead}$(f(\theta),\theta_i,options))$;
\STATE $acc=1$; ($>\epsilon_s$)
\WHILE{$acc>\epsilon_s$}
\STATE $i=i+1$;
\STATE $\theta_f(i) = $\text{Nelder-Mead}$(f(\theta),\theta_f(i-1),options)$;
\STATE $acc = abs(abs(f(\theta_f(i-1))/f(\theta_f(i)))-1)$; 
\ENDWHILE
\STATE \textbf{return} $\theta_f(i)$
\end{algorithmic}
\vspace{0.3cm}

where $\theta_f(i)$ is the solution after optimization $i$, $\theta_i$ the initial solution, $options$ contains the stopping criteria and tolerances of the Nelder-Mead implementation, $abs$ is the absolute value, $\epsilon_s$ is the stopping accuracy required for $acc$, the relative improvement between the last and the current objective value.

We already mention another kind of useful restarts, also suggested e.g. in \cite{LL04}, the global restarts (or multi-starts). These consist simply of running other optimizations from other different initial solutions $\theta_i$. The aim of these global restarts allows to better span the search space which is useful to try to: 1) get around the non-smoothnesses that could have caused other optimizations to stall, for non-smooth objective functions 2) to find better local optima (nearby), for multimodal objective functions 3) get away from suboptimal solutions against the infinite penalization barrier, which is an implicit constraint. By obtaining an increasing number of final solutions having the best objective value found so far (to a given accuracy), we get increasing probabilities that these solutions are locally optimal and in a lesser measure globally optimal.

\subsection{The advantages}

This illustrating algorithm may compete with other methods, e.g. with \cite{GHMO09} using benchmarks of static output feedback optimization in \cite{S11}. It can often find the same objective values or even better ones, possibly in shorter computational times although this is more an exception rather than the rule, the main drawback being the lack of explicit handling of non-smoothnesses as proposed in \cite{AN06,GHMO09}. The main advantage is the great flexibility to handle any objective function $f(\theta)$, only some of which can be written under other frameworks like LMIs. Moreover, even without using gradient information, it can often lead to locally optimal solutions. Also, unlike usual quasi-Newton methods, NM has the ability to explore neighboring valleys with better local optima and likewise this exploring feature may allow NM to overcome non-smoothnesses.

The flexibility is what makes it a candidate of choice for dealing with designs requiring particular structures or properties, such as ensuring the positivity of the solution like desired in this paper, while exploring a non-smooth objective function. More precisely, NM will be particularly efficient for optimization objectives as an alternative to ILMIs, as motivated before and described further as follows.

ILMIs are often efficient for feasibility problems such as finding a stabilizing controller, possibly also ensuring a given performance level, even though they may fail to find a solution even if there exists one (see \cite{I99,EOA97}). However when used for minimization problems, like minimizing a norm of a performance channel, ILMIs have in general no guarantee of convergence to locally optimal solutions (as discussed in \cite{SRSDW11} and noted in \cite{LLS10}). For instance in \cite{LLS10} this stems essentially from the fact that the objective minimization would require a combinatorial optimization problem involving the objective ($\gamma$) and a parameter ($\epsilon$). Considering this important drawback, other methods should be used instead of ILMIs for objective minimizations. ILMIs can be useful when obtaining feasible (e.g. stable) initial solutions is not trivial, which is not the case here (the system is already stable). And then more suited methods should be used for optimizing the objective, like DS methods.

Indeed using these methods can often lead to locally optimal solutions, unlike the algorithm of \cite{LLS10} that produces suboptimal solutions under an a-priori defined level $\gamma$. Although the non-smoothness prevents guarantees of local optimality for direct search methods, it is conjectured at places in the literature that this does not happen often (see ref. in \cite{AN06}). How `often' depends on the density of non-smoothnesses in the objective function and in particular at `partial' optimal solutions, optimal in some but not all directions. This is actually a matter of debate, on which the reader is advised to consult e.g. \cite{AN06,GHMO09} and the references therein. And as noted before, global restarts can be used to try to
reduce the likeliness of getting stuck at a suboptimal solution.

Many optimization problems in systems and control can be quickly implemented to be solved with DS methods. Indeed the objective function is easily built and evaluated with adequate methods, for example the functions from the Control System Toolbox of Matlab. So there is no need to modify or add cumbersome LMIs, which are also especially large with the technique of \cite{LLS10} (see the matrix inequality (16) there). Likewise the unhandy implementation of an ILMI is avoided, moreover often involving parameter(s) for which no automatic a priori choice can be made. In short, the method can be implemented in a brief time even by non-expert users. Also DS methods are useful for simulation-based optimization \cite[Sec. 1.21]{KLT03}, where more intricate objective functions that could hardly be obtained from analytical expressions are computed using simulation results (illustration given in next section). 

An important remark can be made regarding the worst-case computational complexity bounds. In \cite{MCT} we can read that MDS and NM are not competitive with more sophisticated methods such as (quasi-)Newton methods when applied to smooth problems. Indeed it seems clear that gradient expressions and moreover Hessian information should accelerate the search. However a recent result \cite{V11} proves that directional DS methods share the worst-case complexity bound of steepest descent for the unconstrained minimization of a smooth function. This gives further motivation to consider DS methods.

Let us also mention that the formulation of the problem in \cite{LLS10} requires many matrix variables $P_1,P_2,R,U,V,F_i,H_i$, $i =1,...,6$ and three scalars $\alpha,\beta,\gamma$ -with $U,V,\alpha,\beta$ entering non-affinely the matrix inequalities, thus the need of an ILMI- whereas here we only need the original variables $\hat{A},\hat{B},\hat{C},\hat{D}$ which put together lead to the same size as only $R$ in \cite{LLS10}. Using LMIs leads to this typical key problem of inflation of size and number of variables for large systems ($\Sigma$ and $\hat{\Sigma}$). Instead, with general purpose solvers, only the original problem's variables ($\hat{\Sigma}$) are used.
 
LMI problems are solved efficiently with interior point methods converging to the optimum in worst-case polynomial time. This is not guaranteed with the method proposed here, but also not with ILMIs, since the iterative scheme does not preserve the polynomial time complexity. In practice NM deals well with problems with limited number of variables (e.g. $<$ 10-20) but its performances decrease notably for larger number of variables (e.g. $>$ 50-100). As noted before, one should then use the advantage that DS methods can handle problems formulations with the least number of variables. 

\section{Results on example}
\label{NumSim}

The positive system considered is given by \cite{LLS10}:
\[
    x_{k+1} = \begin{bmatrix}
		0.1595 & 0.1890 & 0.2713\\
		0.5091 & 0 & 0\\
		0 & 0.6740 & 0
		\end{bmatrix} x_k +  \begin{bmatrix}
		0.1350 & 0.0128\\
		0.3850\footnotemark[1] & 0.0510\\
		0.1021 & 0.1250
		\end{bmatrix}w_k
\] \footnotetext[1]{This value is erroneously written 0.0128 in \cite{LLS10}.}
\[
    y_k  = 
			\begin{bmatrix}
			0 & 1 & 0\\
			0 & 0 & 1
			\end{bmatrix}x_k + \begin{bmatrix}
			0 & 0.1250\\
			0.1460 & 0
			\end{bmatrix}w_k
\]  
The objective is to estimate $z_k=[1\ 0\ 0]x_k = x_k^1$. For the problem considered, any positive and stable filter can be used as initial solution. Since like in \cite{LLS10} we want to design a first order filter, we can choose the scalar $\hat{A}$ in $[0,1[$ and the other variables in $\mathbb{R}^+$. We then simply use \texttt{rand} of Matlab to generate the random initial solutions (entries chosen uniformly in $[2^{-53}, 1-2^{-53}]$). Performing three optimizations, each from a different random initial solution, gives us the best following solution:

	\[\hat{A} = 0.06978,\ \hat{B} = \left[0.53667\ 2.13004\right], \hat{C} = 0.15218,\ \hat{D} = \left[0.15435\ 0.10931\right]
\]
Instead of the solution in \cite{LLS10}:
	\[\hat{A} = 0.22819,\ \hat{B} = \left[0.00003\ 0.00003\right], \hat{C} = 0.14130,\ \hat{D} = \left[0.17889\ 0.34404\right]
\]
Let it be clear that we do not need the solution of \cite{LLS10} as initial solution, indeed here we simply use independent random feasible solutions and therefore the proposed technique is self-standing. Perhaps for some particular problems ILMIs would be necessary to find initial solutions, but this is not the case here. We nevertheless give as further information what happens when using the solution of \cite{LLS10} as initial solution, at the end of this section.

Note also that since linear systems admit an infinity of equivalent state-space representations, we can restrict the search space to one representation (e.g. by fixing here $\hat{C}=1$). On the one hand, reducing the number of variables reduces the computational time. On the other hand, this removes the possibility for DS methods to explore different state-space coordinates, therefore reducing their ability to find better local optima and getting around non-smoothnesses or away from the border of (implicit) constraints.

\begin{figure}[h]
\begin{center}
\includegraphics[height=9cm]{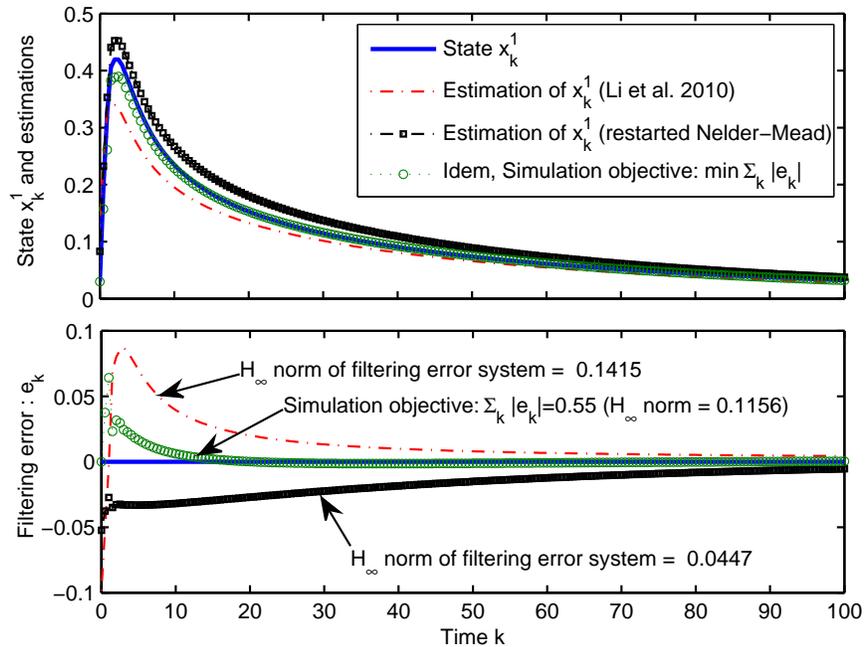} 
\vspace{-0.7cm}   
\caption{State $x_k^1$, estimations and corresponding errors $e_k$.}  
\label{fig1}                                 
\end{center}                                 
\end{figure}

Two of the three solutions found have a $\mathcal{H}_\infty$ performance level around $0.0447$, significantly better than the level $0.1415$ reached in \cite{LLS10}. We present in Fig. \ref{fig1} the same simulation as in \cite{LLS10} of the actual state $x_k^1$ and its estimations, using the following initial condition of the error system: $[0.03, 0.08, 0.10, 0.05]^T$ and exogenous disturbance input:  $w_k = [1/(1+0.25k),e^{-0.02k}]^T$.

The gain of performance level can be seen, with a maximum absolute error about three times smaller with the restarted NM. A third simulation response is also provided, obtained by minimizing an objective specific to this simulation: the sum of the absolute values of the errors $e_k$. The obtained filter is given by: $\hat{A} = 0.011,\ \hat{B} = \left[0.00598\ 1.779\right], \hat{C} = 0.3,\ \hat{D} = \left[0.000046\ 0.0607\right]$. This gives an illustration of simulation-based optimization.

To give further illustration of the performance of the proposed method, we run 100 optimizations each starting from a different random initial solution. We then give the minimum and the average of all 100 obtained objectives values. Also we give the percentage of objectives values that were smaller than 0.1415 and of those smaller than 0.0448. As last indication, we give the average of the 100 computational times required in seconds\footnote{The computer used is a HP Compaq dc7800\copyright, processor Intel Q9300\copyright, 2.5GHz, 3.48Go RAM, software MATLAB 2007b\copyright}. We perform these tests with the proposed technique as well as with other optimizations techniques, each using the same 100 random initial points. The results are given in Table \ref{Result100}.
\begin{table}[h]
\caption{Comparison of objective values and computational times}
\vspace{-0.8cm}
\label{Result100}
\begin{center}
\begin{tabular}{|ll|l|c|c|c|}
\hline
  & Method & Obj.: \texttt{min}/\texttt{mean}  & $<$0.1415 & $<$0.0448 & Time (s)\\
\hline
1) & \texttt{fminsearch} & 0.04470897/0.07 & 96\% & 5\% & 9.57 \\
\hline
2) & \texttt{fminsearch} + restart(s) & 0.04470746/0.04559  & 100\% & 73\% & 44.2 \\
\hline
3) & \texttt{nmsmax} & 0.04470746/0.0514  & 98\% & 38\% & 8.53 \\
\hline
4) & \texttt{nmsmax} + restart(s) & 0.04470746/0.04471  & 100\% & \textbf{100\%} & 41.7 \\
\hline
5) & \texttt{mdsmax} & 0.04471485/0.0479  & 100\% &5\% & 6.81 \\
\hline
6) & \texttt{gridnm} & 0.04470746/0.045 & 100\% & 79\% & 22.8 \\
\hline
7) & \texttt{fminunc} & 0.04473497/0.272 & 56\% & 1\% & 2.41 \\
\hline
8) & \texttt{fmincon} & 0.04470795/0.0627 & 98\% & 38\% & 2.37 \\
\hline
\end{tabular}
\end{center}
\end{table}
\vspace{-0.5cm}

The methods are named after their Matlab implementation file: 1) \texttt{fminsearch}, the Matlab Optimization Toolbox NM implementation 2) The same with local restart(s), as described in the previous section 3) \texttt{nmsmax}, N. Higham NM implementation \cite{MCT} 4) The same with local restart(s) 5) \texttt{mdsmax} N. Higham MDS implementation \cite{MCT} \cite{T89} 6) The grid restrained NM algorithm \texttt{gridnm} \cite{BPT06}, proved convergent on smooth functions 7) \texttt{fminunc} also using the unconstrained formulation (using $10^{100}$ instead of $\infty$ for constraint penalization) 8) \texttt{fmincon}, the only method using the constrained formulation of the problem.\\ 
All accuracies required: `TolF', `TolX', `TolCon' (see Matlab help \cite{Opti} for descriptions), $\epsilon_s$ and the accuracies in \texttt{nmsmax, mdsmax, gridnm} have been set to $10^{-7}$ as well as the tolerance of the objective $||\mathcal{G}||_\infty$ (\texttt{norm(G,inf,tol)}) evaluation. Also the numbers of iterations and functions evaluations were not limited (e.g. `MaxIter' and `MaxFunEvals' set to inf). 

Except for 7), at least 96\% of the trials give a better performance level than 0.1415 from \cite{LLS10}. This is a first illustration of how general purpose optimization solvers might compete with ILMIs. Observe then the results 2), 4) and 6), obtained with improved NM. As can be seen these have high success rate at reaching objectives values very close to the best one found. 4) is especially good since it was 100 \% successful. Considering this, even though $f(\theta)$ is non-smooth here, the solutions found with objective value $<0.0448$ are very likely locally optimal and probably globally as well. This illustrates the efficiency of the method put forward.

The average computational times are under the minute, therefore the methods are not only convenient to encode but also reasonably fast to run. Note that the computational times could be shortened by using gradient expressions like in \cite{AN06,GHMO09} or, as suggested in \cite{AN06}, by interrupting the bisection algorithm computing the $\mathcal{H}_\infty$ norm once it is sure that the solution being evaluated is better or worse than the other solutions. According to one of the authors in \cite{LLS10} the computational time required there was around a few seconds, which is comparable to those required here.

We also run the different alternatives starting from the solution in \cite{LLS10}, which gives the following improved objectives values: 1) 0.1395 2) 0.0565 3) 0.0455 4) 0.0447 5) 0.0466 6) 0.0449 7) 0.1391 8) 0.0454. We see that two methods were unable to improve the solution of \cite{LLS10} to be close to the objective value 0.0447, which gives an indication that this solution is badly located. Indeed its very small $\hat{B}$ matrix, `blocked' near the positivity constraint, almost cancels the effect of the dynamical part of the filter: the solution in \cite{LLS10} can be approximated by a static filter (the $\hat{D}$ matrix) with almost no impact (0.003$\%$) on the objective\footnote{With value actually around 0.1417 and not 0.1415, using the tolerance $10^{-7}$ in \texttt{norm} instead of the default $10^{-2}$}. On this matter we also mention that an important element for the performance of ILMIs is to try to stay away from the borders of the BMI feasible set (using e.g. \cite[Sec. 4.3]{SRSDW11}), otherwise the algorithm may stall early at partial optimal solutions. 

\section{Conclusion}

Considering the heuristical nature of ILMIs, which can be efficient for feasibility problems but have in general no guarantee of convergence to locally optimal solution for minimization objectives nor a bound on the worst-case complexity, many other methods can be used as competitors of such techniques.

The problem of \cite{LLS10}, designing a reduced-order positive filter to estimate the output of a positive system under a given maximum $\mathcal{H}_\infty$ error level $\gamma$, is dealt with in that paper using an ILMI and indeed it can be read there that the convergence of the algorithm is not guaranteed. This motivates us to propose techniques where the $\mathcal{H}_\infty$ error level is minimized. We then consider a less restrictive approach than (I)LMIs and reformulate the problem so as to be solvable by general purpose optimization solvers. In particular, we put forward the NM algorithm improved with local restarts that can behave well even with non-smooth objective functions.

This approach has apparently only advantages compared to that of \cite{LLS10}. It is easy to encode and use even by non-expert users. It is more flexible, i.e. straightforward to modify for example to 1) handle continuous-time systems 2) change the objective into $\mathcal{H}_2$ minimization or multi-objectives 3) take into account complex requirements such as a structure of the solution. The inflation of size of the system -but not necessarily of the filter- will have a smaller impact on the computational time than techniques using LMIs, since the only variables are those of the filter to be designed (no additional variables needed, like Lyapunov matrices). Also the technique often leads to locally optimal solutions, where `often' depends on the objective to be minimized, whereas this should be seldom the case with most ILMIs (see \cite{SRSDW11} and references therein, in particular \cite{HM97}). Finally since it is generally fast for moderate size problems, it can be run multiple times from several feasible initial solutions. And so by getting several times the same best solution one gets an increasing probability that this solution is indeed locally optimal or even globally optimal.

In the end our claim is that most open optimization problems in systems and control should not be handled with conservative LMI formulations, except for providing an initial suboptimal solution. Indeed conservative LMI formulations or ILMIs to deal with BMI problems typically do not have guaranteed convergence to local optima. To get locally optimal solutions one will need to use other approaches, the best two methods being probably \cite{AN06,GHMO09}, which use gradient expressions and take non-smoothnesses into account. Nevertheless the method proposed can be a competitor to these techniques by finding the same or better objectives. Also the set of problems that can be solved by DS methods is broader than those dealt with the methods in \cite{AN06,GHMO09}, thanks to the great flexibility offered by the general objective formulation $f(\theta)$ and without need of gradient expression.

\section*{Acknowledgment}
The authors thanks Ping Li for comments on the manuscript and \'Arp\'ad B\"urmen for providing the grid restrained NM implementation \cite{BPT06}. This research is supported by the Belgian Network DYSCO (Dynamical Systems, Control, and Optimization) funded by the Interuniversity Attraction Poles Programme of the Belgian State, Science Policy Office.

\end{document}